\documentclass[12pt]{amsart}

\newtheorem{thm}{Theorem}

\theoremstyle{definition}

\newtheorem{rem}[thm]{Remark}

\theoremstyle{remark}

\numberwithin{equation}{section}

\def\Q{{\mathbb Q}}
\def\Z{{\mathbb Z}}
\def\C{{\mathbb C}}

\def\Gal{\text{\rm Gal}}

\def\Tor{\text{\rm Tor}}

\begin{document}

\title[Manin-Mumford conjecture]
{On the Manin-Mumford conjecture for abelian varieties
 with a prime of supersingular reduction}

\author[Tetsushi Ito]{Tetsushi Ito}
\address{Department of Mathematics, Faculty of Science, Kyoto University
Kyoto 606-8502, Japan}
\email{tetsushi@math.kyoto-u.ac.jp}

\subjclass{Primary: 14K12; Secondary: 11G10, 14G15}
\date{\today}

\begin{abstract}
We give a short proof of
the \lq\lq prime-to-$p$ version" of the Manin-Mumford conjecture
for an abelian variety over a number field,
when it has supersingular reduction at a prime dividing $p$,
by combining the methods of Bogomolov, Hrushovski, and Pink-Roessler.
Our proof here is quite simple and short,
and neither $p$-adic Hodge theory nor model theory is used.
The observation is that a power of a lift of the Frobenius element
at a supersingular prime acts on the prime-to-$p$ torsion points
via nontrivial homothety.
\end{abstract}

\maketitle

Let $A$ be an abelian variety over a number field $K$.
For an integer $n \geq 1$, let $[n] \colon A \to A$ be
the multiplication-by-$n$ isogeny on $A$,
whose kernel $A[n]$ has $\Gal(\overline{K}/K)$-action.
Let $\Tor(A) := \bigcup_{n \geq 1} A[n]$ be the group of
torsion points on $A$.
For a prime number $p$, let $\Tor_p(A) := \bigcup_{n \geq 1} A[p^n]$
(resp.\ $\Tor^p(A) := \bigcup_{(n,p) = 1} A[n]$)
denote the subgroup of $\Tor(A)$
consisting of elements which have $p$-power order
(resp.\ order prime to $p$).

The Manin-Mumford conjecture states that,
for an irreducible closed subvariety $X$ of $A$,
if $\Tor(A) \cap X $ is Zariski dense in $X$,
then $X$ is a translate of an abelian subvariety of $A$.
For a prime number $p$,
we can also consider the \lq\lq $p$-primary version"
(resp.\ \lq\lq prime-to-$p$ version")
by considering $\Tor_p(A)$ (resp.\ $\Tor^p(A)$)
instead of $\Tor(A)$.
This conjecture was originally proposed by Lang
for curves embedded in their Jacobian varieties (\cite{LangDiv}).
A first partial result was obtained by Bogomolov,
who proved the \lq\lq $p$-primary version" for each $p$
using $p$-adic Hodge theory (\cite{Bogomolov}).
The first full proof was obtained by
Raynaud by rigid analytic methods (\cite{Raynaud1}, \cite{Raynaud2}).
Up to now, several different proofs and generalizations are known
by Hrushovski, Ullmo, Szpiro, Zhang, Pink-Roessler, Roessler
(\cite{Hrushovski}, \cite{Ullmo}, \cite{Zhang},
\cite{PinkRoessler1}, \cite{PinkRoessler2},
\cite{Roessler}.\ see also \cite{Tzermias}).

The aim of this paper is to give a short proof of
the \lq\lq prime-to-$p$ version" of the Manin-Mumford conjecture
for an abelian variety over a number field,
when it has supersingular reduction at a prime dividing $p$,
by combining the methods of Bogomolov, Hrushovski, and Pink-Roessler
(see Remark \ref{RemarkDifference}
for a difference between the proof here and previously known ones).

\begin{thm}
\label{MainTheorem}
Let $A$ be an abelian variety of dimension $g$
over a number field $K$.
Assume that $A$ has supersingular reduction at a prime $v$ dividing $p$.
Let $F$ be the residue field at $v$.
Let $X$ be an irreducible  closed subvariety of $A$ defined over $K$.
If $\Tor^p(A) \cap X $ is Zariski dense in $X$,
then $X$ is a translate of an abelian subvariety of $A$,
that is, $X = A' + a$ where $A'$ is an abelian subvariety of $A$
and $a \in A$.
\end{thm}

\begin{proof}
Let $A_F$ be the reduction of $A$ at $v$,
which is a supersingular abelian variety over $F$.
Let $q$ be the cardinality of $F$, which is a power of $p$.
Let $\sigma \in \Gal(\overline{F}/F)$ be the $q$-th power
Frobenius automorphism of $\overline{F}$.
For a prime number $\ell \neq p$, the characteristic polynomial of $\sigma$
acting on the $\ell$-adic Tate module
$V_{\ell} A_F = \big( \varprojlim_{n} A_F[\ell^n] \big) \otimes_{\Z_{\ell}} \Q_{\ell}$
has coefficients in $\Z$, and is independent of $\ell \neq p$.
Let $\alpha_1,\ldots,\alpha_{2g}$ be the eigenvalues $\sigma$
acting on $V_{\ell} A_F$.
It is known that, for each $i$ and each embedding
$\iota \colon \Q(\alpha_i) \hookrightarrow \C$,
the complex absolute value of $\iota(\alpha_i)$ is $q^{1/2}$
(see e.g.\ \cite{Milne}, \cite{Mumford}).

Since $A_F$ is supersingular,
we know that $\alpha_{i}^2/q$ are algebraic numbers
whose absolute values are 1 at all finite and infinite primes.
Hence $\alpha_{i}^2/q$ are roots of unity
(e.g.\ \cite[Chapter V, \S 1]{LangANT}).
Therefore, there is an integer $m \geq 1$
such that $\alpha_1^{2m} = \cdots = \alpha_{2g}^{2m} = q^{m}$.
Since $\sigma$ acts semisimply on $V_{\ell} A_F$ (\cite{Tate}),
this means that $\sigma^{2m} \cdot x = [q^m](x)$ for all $x \in \Tor^p(A_F)$.
Take a lift $\widetilde{\sigma} \in \Gal(\overline{K}/K)$ of $\sigma$.
Since $\Tor^p(A)$ is naturally isomorphic to $\Tor^p(A_F)$,
we have $\widetilde{\sigma}^{2m} \cdot x = [q^m](x)$
for all $x \in \Tor^p(A)$.
On the other hand, since $X$ is stable
by the action of $\Gal(\overline{K}/K)$,
we have $[q^m](x) \in X$ for all $x \in \Tor^p(A) \cap X$.
Since $\Tor^p(A) \cap X$ is Zariski dense in $X$,
we conclude $[q^m](X) = X$.

Therefore, $X$ is a translate of an abelian subvariety of $A$
by the same argument as in \cite[Th\'eor\`eme 3]{Bogomolov}.
We reproduce the proof here for the completeness of the paper.
After replacing $A$ by its quotient by
$\text{Stab}_A(X) := \{ a \in A \mid X + a = X \}$,
we may assume $\text{Stab}_A(X)$ is trivial.
We shall prove $X$ is a point.
Let $d$ be the dimension of $X$, and put $n := q^m$.
Consider the action of $[n]$ on the cycle class
$[X] \in H^{2(g-d)} \big( A \otimes_K \!\overline{K},\Q_{\ell}(g-d) \big)$
for some $\ell$.
Since $[n] \colon A \to A$ is an \'etale Galois covering of $A$
and the Galois group is $A[n]$ acting on $A$ by translation,
$[n]^{-1}(X) = \bigcup_{a \in A[n]} (X+a)$.
Since $\text{Stab}_A(X)$ is trivial,
we know that $X+a$ are different for all $a \in A[n]$.
Hence we have
$[n]^{\ast}[X] = \sum_{a \in A[n]} [X+a] = n^{2g} [X]$.
Since $[n]^{\ast}$ acts on 
$H^{2(g-d)} \big( A \otimes_K \!\overline{K},\Q_{\ell}(g-d) \big)$
via multiplication by $n^{2(g-d)}$, we have $n^{2(g-d)} = n^{2g}$.
Therefore, we have $d = 0$ and $X$ is a point.
\end{proof}

\begin{rem}
\label{RemarkDifference}
Here we compare the proof above and the proofs previously obtained by
Bogomolov, Hrushovski, and Pink-Roessler.
Bogomolov proved, for each $p$,
there exists an element of $\Gal(\overline{K}/K)$ which
acts on the $p$-adic Tate module
via multiplication by an integer $n>1$
(\cite[Corollaire 2]{Bogomolov}).
The proof of this fact relies on
a study of $p$-adic Galois representations
which have Hodge-Tate decomposition at a prime dividing $p$.
It seems very difficult to determine $\sigma$ and $n$ explicitly.
On the other hand, Hrushovski considered a lift of
the Frobenius element at a prime of good reduction.
However, the Frobenius element does not act on
the Tate module via homothety in general.
Hrushovski overcame this difficulty by using
the model theory of difference fields (\cite{Hrushovski}).
Pink-Roessler, inspired by Hrushovski's work,
obtained a new proof using only classical algebraic geometry
(\cite{PinkRoessler1}, see also \cite{PinkRoessler2}).
Their proof is very short and illuminating
but we need a trick of replacing $A$ by a power of it
(see the proof of \cite[Theorem 1.2]{PinkRoessler1}).
The observation of this paper is that
a power of a lift of the Frobenius element
at a supersingular prime acts on the prime-to-$p$ torsion points
via nontrivial homothety.
\end{rem}

\begin{rem}
It is true but not automatic that
a combination of the \lq\lq $p$-primary version" and
the \lq\lq prime-to-$p$ version" implies
the full Manin-Mumford conjecture.
We need some Galois theoretic results for proving such an implication
(see \cite[\S 8.1]{Raynaud2}.\ see also \cite{Roessler}).
It seems difficult to obtain the full Manin-Mumford conjecture
by simply generalizing the method of this paper
because it seems difficult to control
the action of the Galois group on $\Tor_p(A)$.
Note that, in \cite[\S 3]{PinkRoessler1},
to cover all torsion points,
Pink-Roessler took lifts of the Frobenius elements
at {\it two} primes of good reduction of different characteristics,
and used deep ramification-theoretic results of Serre.
Recently, Roessler succeeded in removing the use of Serre's results
by using a result of Boxall (\cite{Roessler}).
Finally, we note that the idea of finding an element in the Galois group
acting on the Tate module via nontrivial homothety
to prove the Manin-Mumford conjecture is classical, and due to Lang
(\cite[\S 3]{LangDiv},
see also \cite[Conjecture 5.1]{Tzermias}).
\end{rem}

\vspace{0.2in}

\noindent
{\bf Acknowledgments.}
The author would like to thank Seidai Yasuda and Teruyoshi Yoshida
for valuable comments and discussions.
This paper was written during the author's stay at Harvard
in October and November 2004. He would like to thank 
the Department of Mathematics of Harvard University
for its cordial hospitality.
The author was supported by the Japan Society for the
Promotion of Science Research Fellowships for Young
Scientists.

\end{document}